\documentclass[10pt,a4paper]{article}
\title{Spectra of alternating Hilbert operators} 
\author{Nobushige Kurokawa and Hiroyuki Ochiai%
\thanks{2000 Mathematics Subject Classification: 11M06\newline
The second author is supported in part by Grand-in-Aid for 
Scientific Research (A)  No. 19204011.}}
\date{Dedicated to Professor Toshikazu Sunada
for his sixtieth birthday.}
\newcommand{\address}{\noindent
Department of Mathematics, Tokyo Institute of Technology, \\
Oh-okayama, Tokyo 152-8551, Japan. \\
e-mail: kurokawa@math.titech.ac.jp \\
\\
Department of Mathematics, Nagoya University, \\
Chikusa, Nagoya 464-8602, Japan. \\
e-mail:  ochiai@math.nagoya-u.ac.jp 
}
\usepackage{amsmath}
\usepackage{amssymb}
\newtheorem{theorem}{Theorem}
\newtheorem{conjecture}{Conjecture}

\newcommand{\spect}{{\mathrm{Spect}}}
\newcommand{\quant}{{\mathrm{quant}}}
\newcommand{\trace}{{\mathrm{trace}}}
\newcommand{\bF}{{\mathbf F}}
\newcommand{\Z}{{\mathbf Z}}
\newcommand{\C}{{\mathbf C}}
\newcommand{\Det}{{\mathrm{Det}}}

\begin{document}

\maketitle

\begin{abstract}
Spectra of real alternating operators seem to be quite 
interesting from the view point of explaining the Riemann 
Hypothesis for various zeta functions. 
Unfortunately we have not sufficient experiments concerning this theme.
Necessary works would be to supply new examples of spectra 
related to zeros and poles of zeta functions.
A century ago Hilbert (1907) considered a kind of operators 
representing quadratic forms of infinitely many variables.
Demonstrating the calculation of spectra for alternating Hilbert 
operators we hope to present a novel scheme in this paper.
Authors expect this study encourages experts for further studies.
\end{abstract}


\section{Introduction}
In 1907 Hilbert studied the alternating infinite matrix
\[
A=A_\infty = 
 \left(
\frac{1}{m-n} \right)_{m,n\ge1}
\]
as an interesting example relating to quadratic forms of infinitely many
variables;
see 
Weyl \cite{W} (1908).
Hilbert studied the symmetric infinite matrix
\[
 \left(
\frac{1}{m+n-1} \right)_{m,n\ge1}
\]
also.
A few years later
Schur \cite{Sc} (1911) gave a good upper estimate $\pi$ for 
their spectral radius.

\bigskip

We hope to report on our discovery of a periodic nature of the
 spectra of the finite segment
\[
A_N := \left(
\frac{1}{m-n} \right)_{m,n=1,\dots,N}
\]
(zero on the diagonal),
as in the followig conjecture.
It seems that there is no literature concerning this theme
in contrast to the symmetric case 
where we know many studies.

We make the following basic conjecture.
\begin{conjecture}
\label{conj1}
$\spect(A_N)$ is 
asymptotically periodic in the following sense:

For $N$ even,
\[
\spect(A_N) = \{ \pm i \lambda_1^{(N)}, \ldots, \pm i \lambda_{N/2}^{(N)} \}
\]
with $0<\lambda_1^{(N)} < \cdots < \lambda_{N/2}^{(N)}$
satisfying
\[
\lambda_k^{(N)} \sim \frac{2\pi}{N}\left(k-\frac12\right)
\]
as $N \to \infty$.

For $N$ odd,
\[
\spect(A_N) = \{0\} \cup
\{ \pm i \lambda_1^{(N)}, \ldots,\pm i \lambda_{(N-1)/2}^{(N)} \}
\]
with $0<\lambda_1^{(N)} < \cdots < \lambda_{(N-1)/2}^{(N)}$
satisfying
\[
\lambda_k^{(N)} \sim \frac{2\pi}{N} k
\]
as $N \to \infty$.
\end{conjecture}

\bigskip

The next conjecture is a quantum analogue 
($q$-analogue at $q=\zeta_{N}$)
of Conjecture~\ref{conj1}.
\begin{conjecture}\label{conj2}
Let 
\[
A_N^\quant = \left(
\frac{\sin \frac{\pi}{N}}{\sin \frac{\pi}{N}(m-n)}
\right)_{m,n=1,\ldots,N}.
\]
Then
\[
\spect(A_N^\quant) =
\left\{
\begin{array}{ll}
\left\{ \pm i 2 (\sin \frac{\pi}{N}) (k-\frac12) \mid k = 1,\dots, \frac{N}{2} \right \}
& \mbox{for $N$: even} \\
\\
\{0\} \cup \left\{ \pm i 2 (\sin \frac{\pi}{N}) k \mid k = 1,\dots, \frac{N-1}{2} \right \}
& \mbox{for $N$:  odd.}
\end{array}
\right.
\]
\end{conjecture}

Note that Conjecture~\ref{conj2} might imply Conjecture~\ref{conj1}.
For an integer $n \in\Z$,
a $q$-integer $[n]=[n]_q$ is defined to be
$[n] = \displaystyle\frac{q^{n/2} - q^{-n/2}}{q^{1/2}-q^{-1/2}}$.
This is a polynomial in $q^{1/2}$
and the limit $\displaystyle\lim_{q \to 1} [n]_q = n$.
In this sense, $[n]_q$ is understood to be
a deformation (quantization) of an integer $n$.
In our case, $q$ is taken to be the $N$-th primitive root of unity,
$q=\exp(2 \pi i/N)$.
Then $[n] = \dfrac{\sin \frac{\pi n}{N}}{\sin \frac{\pi}{N}}$.
A native deformation of $A_N=\left( \frac{1}{m-n} \right)$
will be $\left( \frac{1}{[m-n]} \right)$,
which is equal to $A_N^\quant$.

We also consider an {\it oscillatory}  version.
For $0 \le \theta \le \frac{\pi}{2}$, let
\[
A_N(\theta) := \left(
\frac{\cos(m-n)\theta}{m-n} \right)_{m,n=1,\dots,N}
\]
(zero on diagonal),
and
\[
B_N(\theta) := \left(
\frac{\sin(m-n)\theta}{m-n} \right)_{m,n=1,\dots,N}
\]
($\theta$ on diagonal).
\begin{conjecture}\label{conj3}
\begin{enumerate}
\item[{\rm(1)}]
The set $\spect(A_N(\theta))$ splits into the following two subsets:
\begin{itemize}
\item
$\left\lfloor \left(1- \frac{\theta}{\pi}\right) N \right\rfloor$ spectra,
which are ``major'' periodic in the interval $[-i(\pi-\theta), i(\pi-\theta)]$.
\item
$\left\lfloor \frac{\theta}{\pi} N \right\rfloor$ spectra,
which are ``minor'' periodic in the interval $[-i\theta, i\theta]$.
\end{itemize}
\item[{\rm(2)}]
The set $\spect(B_N(\theta))$ splits into the following two subsets:
\begin{itemize}
\item
$\left\lfloor \frac{\theta}{\pi} N \right\rfloor$ spectra,
which are ``almost'' $\pi$.
\item
The remaining $\left\lfloor \left(1- \frac{\theta}{\pi}\right) N \right\rfloor$ spectra,
which are ``almost'' zero.

\end{itemize}
\end{enumerate}
\end{conjecture}


\bigskip

First
we notice elementary confirmations of Conjectures~\ref{conj1} -- \ref{conj3}
concerning $\trace(A_N^2)$.

\begin{theorem}\label{theorem1}
$\displaystyle
\lim_{N\to\infty} \frac{\trace(A_N^2)}{N} = -\frac{\pi^2}{3}$.
\end{theorem}

\begin{theorem}\label{theorem2}
$\displaystyle
 \trace((A_N^\quant)^2)
= -\left(\sin\frac\pi{N}\right)^2 \frac{(N-1)N(N+1)}{3}$.
\end{theorem}

\begin{theorem}\label{theorem3}
\begin{itemize}
\item[{\rm(1)}]
$\displaystyle
\lim_{N\to\infty} \frac{\trace(A_N(\theta)^2)}{N} = - \left( \frac{\pi^2}{3} + \theta^2 - \pi \theta \right)$.
\item[{\rm(2)}]
$\displaystyle
\lim_{N\to\infty} \frac{\trace(B_N(\theta)^2)}{N} = \pi\theta$.
\end{itemize}
\end{theorem}
\bigskip
We remark that
these Theorems are compatible with Conjectures.

\noindent\underline{Conjecture~1 $\Rightarrow$ Theorem~1}
\[
\trace((A_N)^2) \sim
- \left(\frac\pi{N}\right)^2\sum_{k=1}^N (N-2k+1)^2
=
-\left(\frac\pi{N}\right)^2 \frac{(N-1)N(N+1)}{3}.
\]

\noindent\underline{Conjecture~2 $\Rightarrow$ Theorem~2}
\[
\trace((A_N^\quant)^2) =
- \left(\sin\frac\pi{N}\right)^2 \sum_{k=1}^N (N-2k+1)^2
=
-\left(\sin\frac\pi{N}\right)^2 \frac{(N-1)N(N+1)}{3}.
\]

\noindent\underline{Conjecture~3(1) $\Rightarrow$ Theorem~3(1)}:
\begin{eqnarray*}
&& \hskip-1cm \trace(A_N(\theta)^2) \\
&\sim& 
-2 \sum_{k=1}^{\lfloor(1-\theta/\pi)N/2\rfloor}
\left(\frac{(\pi-\theta)k}{\lfloor(1-\theta/\pi)N/2\rfloor} \right)^2
-2 \sum_{l=1}^{\lfloor(\theta/\pi)N/2\rfloor}
\left(\frac{\theta}{\lfloor(\theta/\pi)N/2\rfloor} \right)^2 \\
&\sim&
-\frac{(\pi-\theta)^3}{3\pi} N - \frac{\theta^3}{3\pi} N \\
&=& - \left( \frac{\pi^2}{3} + \theta^2 - \pi \theta \right) N.
\end{eqnarray*}
%
\noindent\underline{Conjecture~3(2) $\Rightarrow$ Theorem~3(2)}:
\[
\trace(B_N(\theta)^2) \sim \pi^2 \times \lfloor\frac\theta\pi  N \rfloor \sim \pi \theta N.
\]

We report also the proof of Conjecture~\ref{conj2}:
\begin{theorem}\label{theorem4}
Conjecture~\ref{conj2} is valid.
\end{theorem}
Closing this Introduction,
we briefly explain a possible connection to zeta functions. For 
a prime $p$, let
\[
A^{(p)}=\left(
\frac{1}{ m \log p - n \log p}
\right)_{m,n \ge 1}
= \frac{1}{\log p} \left(
\frac{1}{ m - n}
\right)_{m,n \ge 1}.
\]
Then Conjecture~\ref{conj1}  implies 
\[
\lim_{N \to \infty, N: \mbox{\small odd}} 
\spect(N A_N^{(p)}) = \frac{2 \pi i}{\log p} \Z.
\]
This coincides with the set of poles of the zeta function 
\[
\zeta(s,\bF_p)=(1-p^{-s})^{-1}
\]
of the finite field $\bF_p$. 
This fact may indicate that zeros and 
poles of a zeta function
\[
Z(s) =\prod_{P} (1-N(P)^{-s})^{-1}
\]
are explained by the spectra of the alternating matrix
\[
A=\left(
\left(
\frac{1}{m \log N(P) - n \log N(Q)}
\right)_{m,n \ge 1}
\right)_{P,Q}.
\]

\section{Proof of Theorems~\ref{theorem1} to \ref{theorem3}}
\subsection{Proof of Theorem~\ref{theorem1}}
\begin{eqnarray*}
\trace(A_N^2)
&=&
-2 \sum_{1\le m<n \le N} \frac{1}{(n-m)^2}
\\
&=& 
-2 \sum_{k=1}^{N-1} \frac{N-k}{k^2} \\
&=& -2N \sum_{k=1}^{N-1} \frac{1}{k^2} + 2 \sum_{k=1}^{N-1} \frac1k
\\
&\sim& -N \frac{\pi^2}{3}.
\end{eqnarray*}

\subsection{Proof of Theorem~\ref{theorem2}}
\begin{eqnarray*}
\trace((A_N^\quant)^2)
&=& - \sum_{m,n=1}^N \frac{\sin^2 \frac\pi N}{\sin^2 \frac\pi N(m-n)} \\
&=& -2 \sin^2 \frac\pi N \sum_{k=1}^{N-1} \frac{N-k}{\sin^2 \frac \pi N k} \\
&=& - \sin^2 \frac\pi N \sum_{k=1}^{N-1} 
\left( \frac{N-k}{\sin^2 \frac \pi N k} + \frac{k}{\sin^2 \frac \pi N k} \right) \\
&=& - N \sin^2 \frac\pi N \sum_{k=1}^{N-1} 
 \frac{1}{\sin^2 \frac \pi N k}  \\
&=& -N \sin^2 \frac\pi N \cdot \frac{N^2-1}{3} \\
&=& -\sin^2 \frac\pi N \cdot \frac{(N-1)N(N+1)}{3}. \\
\end{eqnarray*}

\subsection{Proof of Theorem~3(1)}
\begin{eqnarray*}
\trace(A_N(\theta)^2)
&=& -2 \sum_{1\le m < n \le N} \frac{\cos^2(n-m)\theta}{(n-m)^2} \\
&=& -2 \sum_{k=1}^{N-1} \frac{\cos^2(k\theta)}{k^2} (N-k) \\
&=& -2 \left\{
N \sum_{k=1}^{N-1} \frac{\cos^2(k\theta)}{k^2}
- \sum_{k=1}^{N-1} \frac{\cos^2(k\theta)}{k}
\right\}.
\end{eqnarray*}
Here from
\[
\sum_{k=1}^{N-1} \frac{\cos^2(k\theta)}{k} = O(\log N), 
\]
\begin{eqnarray*}
\lim_{N\to\infty} \frac{\trace(A_N(\theta)^2)}{N}
&=& -2 \sum_{k=1}^\infty \frac{\cos^2(k\theta)}{k^2} \\
&=& - \sum_{k=1}^\infty \frac{1+\cos(2k\theta)}{k^2} \\
&=& - \left\{ \frac{\pi^2}{6} + \left( \frac{\pi^2}{6} - \pi \theta + \theta^2 \right) \right\} \\
&=& - \left( \frac{\pi^2}{3} + \theta^2 - \pi \theta \right).
\end{eqnarray*}

\subsection{Proof of Theorem~3(2)}
\begin{eqnarray*}
\trace(B_N(\theta)^2)-\trace(A_N(\theta)^2)
&=&
2 \sum_{1\le m<n\le N}\frac{1}{(m-n)^2} + \theta^2 N \\
&=& 2 \sum_{k=1}^{N-1} \frac{N-k}{k^2} + \theta^2 N.
\end{eqnarray*}
This shows
\[
\lim_{N\to\infty} \frac{\trace(B_N(\theta)^2)-\trace(A_N(\theta)^2)}{N}
= 2 \sum_{k=1}^{\infty} \frac{1}{k^2} + \theta^2
= \frac{\pi^2}{3} + \theta^2.
\]
\section{Proof of Theorem~\ref{theorem4}}
We set $\zeta_{2N}=\exp(\pi i/N)$ 
and $\zeta_{N}=\zeta_{2N}^2=\exp(2 \pi i/N)$.
For $k=1,2,\dots,N$, we have
\begin{eqnarray*}
\sum_{1\le n\le N, n \neq m}
\frac{\sin\frac \pi N}{\sin\frac{\pi}{N}(m-n)} 
\zeta_{2N}^{n(2k-1)}
&=&
\sum_{1\le n\le N, n \neq m}
\frac{(\zeta_{2N}-\zeta_{2N}^{-1})\zeta_{2N}^{n(2k-1)}}{\zeta_{2N}^{m-n}-\zeta_{2N}^{n-m}} 
\\
&=&
(\zeta_{2N}-\zeta_{2N}^{-1}) \zeta_{2N}^{-m}
\sum_{1\le n\le N, n \neq m}
\frac{\zeta_{N}^{kn}}{1-\zeta_{N}^{n-m}} 
\\
&=&
(\zeta_{2N}-\zeta_{2N}^{-1}) \zeta_{2N}^{-m}
\sum_{1\le n\le N-1}
\frac{\zeta_{N}^{k(n+m)}}{1-\zeta_{N}^{n}} 
\\
&=&
(\zeta_{2N}-\zeta_{2N}^{-1}) \zeta_{2N}^{m(2k-1)} 
\sum_{1\le n\le N-1}
\frac{\zeta_{N}^{kn}}{1-\zeta_{N}^{n}} \\
&=&
(\zeta_{2N}-\zeta_{2N}^{-1}) \zeta_{2N}^{m(2k-1)} (k- \frac{N+1}{2}).
\end{eqnarray*}
Note that
the last equality follows from
\[
\sum_{n=1}^{N-1} \frac{1}{1-\zeta_N^n} = \frac{N-1}{2},
\]
and 
\[
\sum_{n=1}^{N-1} \frac{1-\zeta_N^{kn}}{1-\zeta_N^n} = N-k
\qquad \mbox{for } k=1,2,\dots,N-1.
\]
We define an invertible matrix $P_N \in GL(N,\C)$ and 
a diagonal matrix $D_N \in M(N, \C)$ by
\[
P_N = \left(\zeta_{2N}^{m(2n-1)}\right)_{m,n=1,\dots,N},
\quad
D_N=\left( 2i (\sin \frac{\pi}{N})(n-\frac{N+1}{2})   \delta_{mn} \right)_{m,n=1,\dots,N}
\]
then
we have an equality
$A_N^\quant P_N = P_N D_N$.
This proves
\[
\spect(A_N^\quant) =
\left\{ 2 i(\sin \frac{\pi}{N}) (k-\frac{N+1}{2}) \mid k = 1,\dots, N \right \}.
\]

\section{Discussion}
\subsection{Szeg\"o's Theorem and Conjecture~\ref{conj3}}

Szeg\"o \cite{Sz} (1920) proved
the uniform distribution property 
of the eigenvalues  $\lambda_k^{(N)}$ ($k=1,\dots,N$)
of the hermitian Toeplitz operator 
\[
T_N = \left( \begin{array}{cccc}
c_0 & c_{-1} & \dots & c_{-(N-1)} \\
c_1 & c_0 & & \vdots \\
\vdots & & \dots & c_{-1} \\
c_{N-1} & \dots & c_1 & c_0
\end{array}
\right)
\]
in the following form:
\[
\lim_{N\to\infty} \frac{F(\lambda_1^{(N)}) + \cdots + F(\lambda_N^{(N)})}{N}
= \frac{1}{2\pi} \int_0^{2\pi} F(f(x)) dx,
\]
where
$F$ is a suitable test function, 
and $f(x)=\displaystyle\sum_{n} c_n e^{inx}$.
We explain below that Conjecture~\ref{conj3} 
is compatible with Szego's theorem.

\subsubsection{For $A_N(\theta)$}

We check the relation with Conjecture~\ref{conj3}(1)
and Szeg\"o's theorem for
a hermitian Toeplitz  matrix
$i A_N(\theta)=(c_{m-n})_{m,n=1,\dots,N}$.
In this case the entries
$c_n=i \cos n \theta/n$,
and then
\[
f(x) = \displaystyle
\left\{ \begin{array}{ll} 
x & \mbox{ for } 0< x < \theta \\
x-\pi & \mbox{ for } \theta<x <2\pi-\theta \\
x-2\pi & \mbox{ for } 2\pi-\theta<x<2\pi.
\end{array} \right.
\]
This shows
\[
\lim_{N\to\infty} \frac{F(\lambda_1^{(N)}) + \cdots + F(\lambda_N^{(N)})}{N}
= \frac{1}{2\pi} 
\int_{-\pi}^\pi F(x) m(x) dx,
\]
where
\[
m(x) = \displaystyle
\left\{ \begin{array}{ll} 
2 & \mbox{ for } |x| < \theta \\
1 & \mbox{ for } \theta<|x| <\pi-\theta \\
0 & \mbox{ for }  \pi-\theta<|x|<\pi.
\end{array} \right.
\]
This is compatible with Conjecture~\ref{conj3}(1).

\subsubsection{For $B_N(\theta)$}

In this case 
the hermitian Toeplitz matrix is
$B_N(\theta)= (c_{m-n})_{m,n=1,\dots,N}$ with $c_n = \sin n \theta/ n$.
Using the formula
\begin{eqnarray}
\sum_{n=1}^\infty \frac{\sin nx}{n} = \frac{\pi-x}{2}
\quad \mbox{for } 0 < x < 2 \pi,
\end{eqnarray}
we have
\[
f(x) = \displaystyle
\left\{ \begin{array}{ll} 
\pi & \mbox{ for } 0 < x < \theta \\
0 & \mbox{ for } \theta<x<2\pi-\theta \\
\pi & \mbox{ for } 2\pi- \theta<x < 2\pi. \\
\end{array} \right.
\]
This shows
\[
\lim_{N\to\infty} \frac{F(\lambda_1^{(N)}) + \cdots + F(\lambda_N^{(N)})}{N}
= \frac{\theta}{\pi} F(\pi) + \frac{\pi-\theta}{\pi} F(0).
\]
This is compatible with Conjecture~\ref{conj3}(2).

\subsection{Expected relation to zeta functions}
We explain a possible way  to reach
\[
A=\left(
\left(
\frac{1}{m \log N(P) - n \log N(Q)}
\right)_{m,n \ge 1}
\right)_{P,Q}
\]
from a zeta function
\[
Z(s) =\prod_{P} (1-N(P)^{-s})^{-1}.
\]
As a typical example of a zeta function 
we take up the Riemann zeta function $\zeta(s)$.
Suppose that we have a determinant expression
\[
\hat\zeta(s) \cong \frac{\Det (A-(s-\frac12))}{s(s-1)}
\]
for the completed Riemann zeta function
\[
\hat\zeta(s) = \zeta(s) \pi^{-\frac{s}{2}} \Gamma\left(\frac s2 \right)
\]
with a real alternating matrix $A$.
Then we will have a proof of Riemann Hypothesis
as suggested by Hilbert and Polya around 1915.

Now, on the other hand we have
\[
\hat\zeta(s) = \frac{e^{as+b}}{s(s-1)} \prod_{\rho} 
\left( 1-\frac s\rho \right) e^{s/\rho},
\]
where $\rho$ runs over essential zeros.
The logarithmic derivation gives
\[
\frac{\zeta'}{\zeta}(s) - \frac 12 \log \pi 
+ \frac12 \frac{\Gamma'}{\Gamma}\left( \frac s2 \right)
= a - \frac 1s - \frac 1{s-1} + \sum_{\rho} 
\left( \frac1{s-\rho} + \frac1\rho \right).
\]
Hence we have
\[
\left( \frac{\zeta'}{\zeta} \right)'(s)
+ \frac14 \left(\frac{\Gamma'}{\Gamma}\right)' \left( \frac s2 \right)
= \frac1{s^2} + \frac1{(s-1)^2} - \sum_\rho \frac1{(s-\rho)^2}.
\]
In particular,
\[
\sum_\rho \frac1{(\rho-\frac12)^2}
= - \left( \frac{\zeta'}{\zeta} \right)'\left(\frac12\right)
- \frac14 \left(\frac{\Gamma'}{\Gamma}\right)' \left( \frac14 \right)
+8.
\]
Hence, elementary calculation shows that
\begin{eqnarray*}
\left(\frac{\Gamma'}{\Gamma}\right)'\left(\frac14\right)
&=& \sum_{n=0}^\infty \frac{1}{(n+\frac14)^2} \\
&=& 16 \sum_{n=0}^\infty \frac{1}{(4n+1)^2} \\
&=& 16 \cdot \frac12 \left\{
(1-2^{-s}) \zeta(2) + L(2,\chi_{-4}) \right\} \\
&=& \pi^2 + 8 G,
\end{eqnarray*}
where
\[
G=L(2,\chi_{-4}) = \sum_{n=0}^\infty \frac{(-1)^n}{(2n+1)^2}
\]
is the Catalan constant.
Thus 
\[
\sum_{\rho} \frac{1}{(\rho-\frac12)^2}
= - \left( \frac{\zeta'}{\zeta} \right)' \left( \frac12 \right)
- \frac{\pi^2}{4} - 2G + 8.
\]
Here it might be suggestive to write this as
\[
\sum_{\rho} \frac{1}{(\rho-\frac12)^2}
= - \mbox{``$\displaystyle
\sum_{p,m \ge 1} ( \log p)^2 m p^{-\frac m2}
 $''}
- \frac{\pi^2}{4} - 2G + 8
\]
since we see
\[
\left( \frac{\zeta'}{\zeta} \right)' \left(  s \right)
= \sum_{p,m} (\log p)^2 m p^{-ms}
\]
at least for $\Re(s) > 1$.
Of course we have similarly the formula
\[
\sum_{\rho} \frac{1}{(\rho-\frac12)^{2k}}
= - \mbox{``$\displaystyle
\frac{1}{(2k-1)!} \sum_{p,m} ( \log p)^{2k} m^{2k-1} p^{-\frac m2}
 $''}
+ \alpha_k
\]
with
\[
\alpha_k = -\frac12 (2^{2k}-1) \zeta(2k) - 2^{2k-1} L(2k,\chi_{-4}) + 2^{2k+1}
\]
for $k=1,2,3,\dots$.
On the other hand,
this should be equal to the trace
\[
\mbox{Trace}(A^{-2k})
\]
because of our determinant expression.
Thus we should have
\[
\mbox{Trace}(A^{-2k}) = 
\mbox{``$\displaystyle
\frac{1}{(2k-1)!} \sum_{p,m} ( \log p)^{2k} m^{2k-1} p^{-\frac m2}
 $''}
+ \alpha_k
\]
for $k=1,2,3,\dots$.
[And moreover we see easily that this trace formula is essentially equivalent to
the determinant expression supposed first.]

Thus, we reach to the basis problem:
determine the index set (basic set) $X$ representing
\[
A=\left(a(i,j) \right)_{i,j\in X}.
\]
Since, formally we may write
\[
A^{-2k} = \left( a_{-2k}(i,j) \right)_{i,j \in X}
\]
and
\[
\mbox{Trace}(A^{-2k}) = \sum_{i \in X} a_{-2k}(i,i),
\]
we would have
\[
\sum_{i \in X} a_{-2k}(i,i) = 
\mbox{``$\displaystyle
\frac{1}{(2k-1)!} \sum_{p^m} ( \log p)^{2k} m^{2k-1} p^{-\frac m2}
 $''}
+ \alpha_k.
\]
Hence,
we reach to an obvious suggestion
\[
X = \{ p^m \mid p: \mbox{ prime}, m \ge 1\},
\]
and a simple
\[
A= \left( \frac{1}{\log(p^m) - \log (q^n)} \right)_{p^m,q^n \in X}
\]
in considering the case of $\zeta(s,\bF_p)$.
Similarly we might expect that
\[
A= \left( \frac{1}{\log (N(P)^m) - \log (N(Q)^n)} \right)_{P^m,Q^n}
\]
would explain central essential zeros and poles of 
\[
Z(s) = \prod_{P} (1-N(P)^{-s})^{-1}
\]
in the sense
\[
Z(s) \cong \Det \left(A-\left(s-\frac{\dim}{2}\right)\right)^{(-1)^{\dim+1}},
\]
where $\dim$ is the Kronecker dimension,
for example:
\[
\dim = \left \{
\begin{array}{ll}
1 & \mbox{for } \zeta(s) \\
0 & \mbox{for } \zeta(s,\bF_p).
\end{array}
\right.
\]

\subsection{Symmetric case}

In the case of the symmetric Hilbert matrix
\[
S_N = \left(\frac{1}{m+n-1}\right)_{m,n=1,\ldots,N}
\]
the spectra $\spect(S_N)$ is contained in $(0, \pi)$ 
and is not periodic. 
But curiously enough a quantization
\[
S_N^\quant = 
\left(\frac{\sin\frac{\pi}{N}}{\sin \frac{\pi}{N}(m+n-1)} \right)_{m,n=1,\dots,N}
\]
(zero on the anti-diagonal) has the periodic spectra.
In fact we have 
\begin{theorem} 
\[
\spect(S_N^\quant)=\{
(N-1)\sin\frac{\pi}{N},
(N-3)\sin\frac{\pi}{N},\ldots,
(3-N)\sin\frac{\pi}{N},
(1-N)\sin\frac{\pi}{N}
\} .
\]
\end{theorem}
Proof.
Let $C_N^\quant=\left(\frac{1}{\sin \frac{\pi}{N}(m+n-1)} \right)_{m,n=1,\dots,N}$
(zero on the diagonal).
We will prove
\begin{eqnarray*}
\spect(C_N^\quant) &=&
\left\{ n \in \Z \mid n - N \mbox{ is odd, and} \left|n\right| < N \right\} \\
&=& \left\{ N-1,N-3,\dots,3-N,1-N \right \},
\end{eqnarray*}
which is the set of weights of the irreducible $N$-dimensional
representation of $SL_2$.

We define a matrix 
\[
Q_N = \left( \cos(\frac\pi{2N} (2m-1)(2n-1)-\frac\pi4 ) \right)_{m,n=1,\dots,N}
\]
and 
a diagonal matrix 
\[
D'_N=\left((N+1-2m)   \delta_{mn} \right)_{m,n=1,\dots,N}.
\]
We claim that $Q_N$ is invertible and
we have an equality
$C_N^\quant Q_N = Q_N D'_N$.

We use
\begin{eqnarray*}
\sum_{n=1}^{N-1} \frac{\cos(\frac\pi N n(2k-1))}{\sin \frac\pi N n} &=& 0, \\
\sum_{n=1}^{N-1} \frac{\sin(\frac\pi N n(2k-1))}{\sin \frac\pi N n} &=& N+1-2k
\end{eqnarray*}
for $k=1,2,\dots,N$.

Lastly, we notice that $\det Q_N = \pm (N/2)^{N/2}$,
where the sign $-$ for $N+1 \in 4\Z$ and the sign $+$ for otherwise.


\address
\end{document}